\documentclass[13pt]{article}
\usepackage[cp1251]{inputenc}
\usepackage[english]{babel}
\usepackage{amssymb,amsmath,amsthm,indentfirst}
\usepackage{multicol}
\usepackage{graphicx}

\begin{document}

\begin{center}
\LARGE{Exact Bounded Boundary Controllability to Rest for the Two-Dimensional Wave Equation}
\end{center}

\vskip 0.5cm
\begin{center}
{\large Igor Romanov}
\footnote{National Research University Higher School of Economics,\\
20 Myasnitskaya Ulitsa, Moscow 101000, Russia}
{\footnote{
E-mail: \it{ivromm1@gmail.com}}}
{\large Alexey Shamaev}
\footnote{Institute for Problems in Mechanics RAS,\\
101 Prosp. Vernadskogo, Block 1, Moscow 119526, Russia}
\footnote{Lomonosov Moscow State University,\\
GSP-1, Leninskie Gory, Moscow 119991, Russia}

\end{center}




\begin{abstract}
The problem of the exact bounded control of oscillations of the two-dimensional membrane is considered. Control force is applied to the boundary of the membrane, which is located in a domain on a plane. The goal of the control is to drive the system to rest in a finite time.
\end{abstract}


\section{Introduction}

The problem of exact boundary controllability of oscillations of a plane membrane is considered. Control force has a restriction on its absolute value.
We will prove that the plane membrane can be driven to rest in a finite time. Exact mathematical definitions will be provided. It should be noted that the given method for the proof in this article can be used in the case of any other dimension, but here the two-dimensional case is provided for clear and simple presentation.

The problem of full stabilization in a finite time in case of the distributed
control is described in the monograph \cite{Chernousko}. This reference also contains the upper estimate for an
optimal control time.

Previously the question of the control of oscillations of a plane membrane by means of boundary forces is considered by many authors (i. g. overviews of D. L. Russell \cite{Russell} and J. Lions
\cite{Lions}, as well as the literature provided there).
The monograph \cite{Butk}
describes the task of stabilizing the oscillations of a restricted string by means of the boundary control, and
proves that vibrations of the string can be driven to rest in a finite time
under the condition of restriction imposed on an absolute value of the control function, and an estimate is
provided for the time that is necessary for full rest.
In monograph \cite{Lions1} problems of the optimal control of systems with distributed parameters are studied and  optimal conditions are stated, which are similar to conditions for systems with a finite number of freedom's degrees. Although this method does not provide a constructive technique for finding an optimal control in many cases. In synoptic article \cite{Lions} the problem of exact zero-controllability of a membrane is considered, the existence of the boundary control is proven and the time estimate is given which is required for driving to rest. Here authors, while studying the problem in various formulations, often reject the requirement of optimality of the control and solve only the problem of controllability, which is much easier. What is more, problems with restrictions of the force's absolute value are not considered, explicit forms for control functions are not found, and only theorems of existence are proven.

The statement of the problem in the article essentially differs from the one in \cite{Russell}
and \cite{Lions}, because the value of control force on the boundary has to satisfy the condition: $|u(t,x)|\leqslant\varepsilon$. Note, here the aim is to find not an optimal control, but the admissible
(satisfying initial restrictions) control.

\section{The statement of the problem}

Let us consider the initial-boundary value problem for the two-dimensional wave equation:
\begin{equation}
\label{1}
w_{tt}(t,x)-\Delta w(t,x)=0,\quad (t,x)\in Q_{T}=(0,T)\times\Omega,
\end{equation}
\begin{equation}
\label{2}
w|_{t=0}=\varphi(x),\quad w_t|_{t=0}=\psi(x),\quad x\in\Omega,
\end{equation}
\begin{equation}
\label{3}
\frac{\partial w}{\partial{\nu}}=u(t,x),\quad (t,x)\in\Sigma,
\end{equation}
where $\Omega\subset R^2$ is a bounded, star-shaped relatively some ball domain with an infinitely smooth boundary, $\nu$ --- the outer normal to the boundary of the domain $\Omega$, $\Sigma$ is a lateral surface of a cylinder $Q_T$. Initial data $\varphi(x)$ and $\psi(x)$ are given and will be
chosen in suitable Hilbert spaces, $u(t,x)$ is a control function defined on the boundary
$\Gamma=\partial\Omega$.

Let $\varepsilon>0$ be an given arbitrary number. Let us impose the constraint on the control function:
\begin{equation}
\label{4}
|u(t,x)|\leqslant\varepsilon,
\end{equation}

The problem is to construct a control $u(t,x)$ satisfying inequality (\ref{4})
such that the corresponding solution $w(t,x)$ to the initial-boundary value problem (\ref{1})---(\ref{3})
and its derivative with respect to $t$ become $(C,0)$ at some time $T$,
i.e.
\begin{equation}
\label{5}
w(T,x)=C,\:\: w_t(T,x)=0,
\end{equation}
for all $x\in\Omega$. In this case $C$ is some constant. If we obtained a control $u(t,x)$ such that
conditions (\ref{5}) are achieved then
the system (\ref{1})---(\ref{3}) is called \emph{controllable to rest}.

The following theorem is the main result of this article.

\textbf{Theorem 1.}
Let $\varphi(x)\in H^6(\Omega)$ and $\psi(x)\in H^5(\Omega)$ such that
$$
\frac{\partial\varphi(x)}{\partial\nu}=\Delta\varphi(x)=
\frac{\partial\Delta\varphi(x)}{\partial\nu}=\Delta^2\varphi(x)=
\frac{\partial\Delta^2\varphi(x)}{\partial\nu}=0,\:\: x\in\Gamma,
$$
\begin{equation}
\label{6}
\psi(x)=\frac{\partial\psi(x)}{\partial\nu}=
\Delta\psi(x)=\frac{\partial\Delta\psi(x)}{\partial\nu}=
\Delta^2\psi(x)=0,\:\: x\in\Gamma.
 \end{equation}
Then the system (\ref{1})---(\ref{3}) is controllable to rest.

The proof of Theorem 1 consists of two steps. The first step stabilizes the considered solution and its first derivative with respect to $t$ in a small vicinity of equilibrium $(C,0)$ in the norm of
$C^4(\overline{\Omega})\times C^3(\overline{\Omega})$, and the second step allows
to drive to rest the system in this small vicinity.

\section{The first step of the control}

Here we state the task to stabilize a pair $(w(t,x),w_t(t,x))$ in an arbitrarily small vicinity of $(C,0)$
in the norm of the space $\mathcal{C}^{4,3}(\overline{\Omega})=C^4(\overline{\Omega})\times C^3(\overline{\Omega})$
where $w(t,x)$ is the solution to the system (\ref{1})---(\ref{3}) and $w_t(t,x)$ is its first derivative with respect to $t$.
What is more the control function should satisfy the restriction (\ref{4}).

At the first step we state the problem to stabilize solution of (\ref{1})---(\ref{3}) and its
first derivative by $t$ to the small enough vicinity of $(C,0)$ by the norm of
$\mathcal{H}^6(\Omega)=H^6(\Omega)\times H^5(\Omega)$.
In this case we have the restriction (\ref{4}) for the control function.

For this purpose we use \cite{Quinn} and \cite{Lagnese}. In these articles authors
consider a friction on $\Gamma$ which is defined by $w_t(t,x)$. More exactly they consider
the initial-boundary value problem (\ref{1})---(\ref{2}) with a new boundary condition:
\begin{equation}
\label{6}
\frac{\partial w(t,x)}{\partial{\nu}}=-k\frac{\partial w(t,x)}{\partial t},\quad x\in\Gamma,
\end{equation}
where $k>0$ is a friction coefficient. Let us illuminate shortly
the questions of solvability of this problem.

Let us denote
$$
H=L_2(\Omega),\quad V=H^1(\Omega),
$$

Let us define in the space $V\times H$ an unbounded operator
$$
\mathfrak{A}=
\left(
  \begin{array}{cc}
    0 & I \\
    \Delta & 0 \\
  \end{array}
\right)
$$
with the domain
$$
D(\mathfrak{A})=\{(w_1,w_2)\in H^2(\Omega)\times H^1(\Omega):
\:\frac{\partial w_1}{\partial\nu}=-kw_2,
\:x\in\Gamma\}.
$$

It is a well known fact that the norm in the space $D(\mathfrak{A})$
can be represent in the following form:
\begin{equation}
\label{10.4}
\|(w_1,w_2)\|_{D(\mathfrak{A})}=\|(w_1,w_2)\|_{V\times H}+
\|\mathfrak{A}(w_1,w_2)\|_{V\times H}.
\end{equation}

Let us consider the following system of differential equations:
\begin{equation}
\label{10.2}
\bar{w}_t=\mathfrak{A}\bar{w},
\end{equation}
where $\bar{w}=(w_1,w_2)$.

It is known (see \cite{Quinn}, \cite{Lagnese}) that an operator $\mathfrak{A}$ is a generator of
strongly continuous semigroup of linear bounded operators.

It is a well known fact that if initial data
$(\varphi,\psi)$ is an element of $D(\mathfrak{A}^k)$, $k=0,1,2,...$,
then we have:
$$
(w_1(t),w_2(t))\in
C\left([0,T];D(\mathfrak{A}^k)\right).
$$
We note that in our case we have $(\varphi,\psi)$ as an element of $D(\mathfrak{A}^5)$.

Let $(\varphi,\psi)\in V\times H$. It is proved
(see \cite{Quinn}, \cite{Lagnese})
that for the energy of the system we have:
\begin{equation}
\label{Energy}
E(t)\rightarrow 0,\:\: t\rightarrow+\infty.
\end{equation}
where
$$
E(t)=\int\limits_{\Omega}\left\{w^2_{1,x_1}(t,x)+w^2_{1,x_2}(t,x)+w^2_{2}(t,x)\right\}dx
$$
is an energy of the system.

We introduce:
$$
C=\frac{1}{|\Gamma|}\int\limits_{\Gamma}\varphi(x)d\Gamma+
\frac{1}{k|\Gamma|}\int\limits_{\Omega}\psi(x)dx,
$$
where $|\Gamma|$ is a length of $\Gamma$.
Let $w(t,x)=v(t,x)+C$ and consider a new initial-boundary value problem for $v(t,x)$
(analogous to (\ref{1}), (\ref{2}), (\ref{6})):
\begin{equation}
\label{v1}
\bar{v}_t=\mathfrak{A}\bar{v},
\end{equation}
\begin{equation}
\label{v2}
(v_1,v_2)|_{t=0}=(\varphi(x)-C,\psi(x)),\quad x\in\Omega,
\end{equation}
where $\bar{v}=(v_1,v_2)$. Obviously in this case $v_1=v$ and $v_2=v_t$.

Using Friedrichs' (Poincare) inequality (see \cite{Mikhlin}) we have
$$
\int\limits_{\Omega}v^2(t,x)dx\leqslant C_3\left\{\int\limits_{\Omega}\left(\left(
\frac{\partial v}{\partial x_1}\right)^2+\left(\frac{\partial v}{\partial x_2}\right)^2
\right)dx+\left(\int\limits_{\Gamma}v(t,x)d\Gamma\right)^2\right\}.
$$
$$
\int\limits_{\Gamma}v(t,x)d\Gamma=-\frac{1}{k}\int\limits_{0}^{t}\int\limits_{\Gamma}
\frac{\partial v}{\partial\nu}d\Gamma dt+\int\limits_{\Gamma}\varphi(x)d\Gamma
-|\Gamma|C=
$$
$$
=-\frac{1}{k}\int\limits_{0}^{t}\int\limits_{\Omega}
\Delta v(t,x)dx dt+\int\limits_{\Gamma}\varphi(x)d\Gamma
-|\Gamma|C=
$$
$$
=-\frac{1}{k}\int\limits_{0}^{t}\int\limits_{\Omega}
v_{tt}(t,x)dx dt+\int\limits_{\Gamma}\varphi(x)d\Gamma
-|\Gamma|C=
$$
$$
=-\frac{1}{k}\int\limits_{\Omega}
v_t(t,x)dx+\frac{1}{k}\int\limits_{\Omega}\psi(x)dx+
\int\limits_{\Gamma}\varphi(x)d\Gamma
-|\Gamma|C=-\frac{1}{k}\int\limits_{\Omega}v_t(t,x)dx.
$$
From the last estimations we obtain
\begin{equation}
\label{L2}
\|w(t,\cdot)-C\|_{L_2(\Omega)}\rightarrow 0,\quad t\rightarrow+\infty.
\end{equation}

Let $(\varphi-C,\psi)$ be an element of $D(\mathfrak{A})$ and
$(v_1(t),v_2(t))$ is a corresponding (to these initial data) solution.
We consider now the following Cauchy problem:
$$
\frac{d}{dt}\mathfrak{A}\bar{v}(t)=\mathfrak{A}^2\bar{v}(t),\:\:
\mathfrak{A}\bar{v}(0)=\mathfrak{A}(\varphi-C,\psi).
$$
We note that
\begin{equation}
\label{10.6}
\mathfrak{A}(v_1(t),v_2(t))=(v_2(t),\Delta v_1(t)).
\end{equation}

Then from (\ref{Energy}) and (\ref{10.6}) we obtain
\begin{equation}
\label{10.7}
\int\limits_{\Omega}\left\{v^2_{2,x_1}(t)+v^2_{2,x_2}(t)+(\Delta v_1(t))^2\right\}dx
\rightarrow 0,\:\: t\rightarrow+\infty.
\end{equation}

Combining (\ref{10.4}), (\ref{Energy}), (\ref{L2}) and (\ref{10.7}),
we have:
\begin{equation}
\label{Energy2}
\|(v_1(t),v_2(t))\|_{D(\mathfrak{A})}\rightarrow 0,\:\: t\rightarrow+\infty.
\end{equation}

Let initial condition is an element of $D(\mathfrak{A})$ then for the
corresponding solution we can obtain (using the theory of elliptic boundary
value problems (see, for example, \cite{Agranovich} or \cite{Madgenes})) the following estimate:
\begin{equation}
\label{10.8}
\|v_1(t)\|_{H^2(\Omega)}\leqslant N_1\left(\|\Delta v_1(t)\|_{L_2(\Omega)}+
k\|v_2(t)\|_{H^{\frac{1}{2}}(\Gamma_1)}+\|v_1(t)\|_{L_2(\Omega)}\right),
\end{equation}
where $N_1$ does not depend on $(v_1,v_2)$.
Using (\ref{Energy}), (\ref{10.7}) and the last estimate
one can easily prove that $v_1(t)$ tends to zero when $t\rightarrow+\infty$
in the norm of $H^2(\Omega)$.

Consider the space $D(\mathfrak{A}^2)$. Using the theory of
elliptic boundary value problems we can describe this space effectively:
$$
D(\mathfrak{A}^2)=\{(v_1,v_2)\in H^3(\Omega)\times H^2(\Omega):\:
\frac{\partial v_1}{\partial\nu}=-kw_2,\:\frac{\partial v_2}{\partial\nu}=
-k\Delta v_1,\:x\in\Gamma\}.
$$

Let $(v_1(t),v_2(t))$ be
the solution to (\ref{1}), (\ref{2}), (\ref{6}) then $(v_1(t),v_2(t))$ is an
element of $C\left([0,T];D(\mathfrak{A}^2)\right)$.
We have
\begin{equation}
\label{D(u^2)}
\mathfrak{A}^2(v_1,v_2)=(\Delta v_1,\Delta v_2).
\end{equation}
It follows from (\ref{D(u^2)}) that
\begin{equation}
\label{10.11}
\int\limits_{\Omega}\left\{(\Delta v_{1,x_1}(t))^2+(\Delta v_{1,x_2}(t))^2
+(\Delta v_2(t))^2\right\}dx\rightarrow 0,\:\: t\rightarrow+\infty.
\end{equation}

Using (\ref{Energy2}) and (\ref{10.11}) we obtain:
\begin{equation}
\label{Energy3}
\|(v_1(t),v_2(t))\|_{D(\mathfrak{A}^2)}\rightarrow 0,\:\: t\rightarrow+\infty.
\end{equation}

The theory of elliptic boundary value problems gives us the following
estimates:
\begin{equation}
\label{10.9}
\|v_1(t)\|_{H^3(\Omega)}\leqslant N_2\left(\|\Delta v_1(t)\|_{H^1(\Omega)}+
k\|v_2(t)\|_{H^{\frac{3}{2}}(\Gamma_1)}+\|v_1(t)\|_{L_2(\Omega)}\right).
\end{equation}

\begin{equation}
\label{10.10}
\|v_2(t)\|_{H^2(\Omega)}\leqslant N_3\left(\|\Delta v_2(t)\|_{L_2(\Omega)}+
k\|\Delta v_1(t)\|_{H^{\frac{1}{2}}(\Gamma_1)}+\|v_2(t)\|_{L_2(\Omega)}\right).
\end{equation}
Using the last estimates
one can easily prove that $v_1(t)$ tends to zero when $t\rightarrow+\infty$
in the norm of $H^3(\Omega)$.

Let us have a look at one more step in detail. Consider the space $D(\mathfrak{A}^3)$.
We have: $\mathfrak{A}^3(v_1,v_2)=(\Delta v_2,\Delta^2 v_1)$. Hence we obtain two equations: $\Delta v_2=f_1$, $\Delta^2 v_1=f_2$, where $(f_1,f_2)\in H^1\times L_2$,
and three boundary conditions at $\Gamma$:
\begin{equation}
\label{BoundaryCon}
a)\: \frac{\partial v_1}{\partial\nu}=-kv_2,\:
b)\: \frac{\partial v_2}{\partial\nu}=-k\Delta v_1,\:
c)\: \frac{\partial\Delta v_1}{\partial\nu}=-k\Delta v_2.
\end{equation}
Let us make a substitution $\Delta v_1=h$, then the equation $\Delta^2 v_1=f_2$ with the boundary condition (c)
has the form $\Delta h=f_2$, $\frac{\partial h}{\partial\nu}=-k\Delta v_2$. Hence the following estimation takes place
$$
\|h\|_{H^2(\Omega)}\leqslant N_4\left(\|\Delta h\|_{L_2(\Omega)}+k\|\Delta v_2\|_{H^{\frac{1}{2}}(\Gamma)}+
\|h\|_{L_2(\Omega)}\right).
$$
Then $\Delta v_1\in H^{\frac{3}{2}}(\Gamma)$ at the boundary of a domain.
So from the equation $\Delta v_2=f_1$ and the boundary condition (b)
the following estimation is derived:
\begin{equation}
\label{10.11}
\|v_2\|_{H^3(\Omega)}\leqslant N_5\left(\|\Delta v_2\|_{H^1(\Omega)}+k\|\Delta v_1\|_{H^{\frac{3}{2}}(\Gamma)}+
\|v_2\|_{L_2(\Omega)}\right).
\end{equation}

Then we get the equation $\Delta v_1=f_3\in H^2(\Omega)$ with the boundary condition (a).
Using the previous estimation, we obtain:
\begin{equation}
\label{10.12}
\|v_1\|_{H^4(\Omega)}\leqslant N_6\left(\|\Delta^2 v_1\|_{L_2(\Omega)}+k\|v_2\|_{H^{\frac{5}{2}}(\Gamma)}+
k\|\Delta v_2\|_{H^{\frac{3}{2}}(\Gamma)}+\|v_1\|_{L_2(\Omega)}\right).
\end{equation}

Continuing in the analogous way we can prove that
$\|v_1(t)\|_{H^6(\Omega)}$ and $\|v_2(t)\|_{H^5(\Omega)}$ tend to zero when $t\rightarrow+\infty$.
It means that
$$
\|w(t)-C\|_{H^6(\Omega)},\:\|w_t(t)\|_{H^5(\Omega)}\rightarrow 0\quad
t\rightarrow+\infty.
$$

We solve the problem (\ref{1}), (\ref{2}), (\ref{6}) with the given initial conditions, then this solution is substituted to the \emph{only} right part of the equality (\ref{6}), and we obtain the boundary condition (\ref{3})
for the initial-boundary value problem (\ref{1})---(\ref{3}). In other words, we make the control function of the problem (\ref{1})---(\ref{3}) be equal to
$$
u^{(1)}(t,x)=-k\frac{\partial w_0(t,x)}{\partial t},
$$
where $w_0$  is a solution to the problem (\ref{1}), (\ref{2}), (\ref{6}).

Therefore it is proved (here we use Sobolev embedding theorem) that controlling for a long time, we can make the values
$$
\|w(T_1,\cdot)\|_{C^4(\overline{\Omega})},\:\:
\|w_t(T_1,\cdot)\|_{C^3(\overline{\Omega})}
$$
arbitrarily close to $(C,0)$ at some time $t=T_1$.

Now let us show that the boundary control function $u(t,x)$ can be sufficiently small,
i.e. we may satisfy the restriction (\ref{4}). It is known that
$$
\max\limits_{t\in[0,+\infty)}E(t)=E(0)=
\int\limits_{\Omega}\left(\varphi^2_{x_1}(x)+\varphi^2_{x_2}(x)+\psi^2(x)\right)dx.
$$
Then using S. L. Sobolev theorems of injections and (\ref{10.10}) we have:
$$
\|w_2(t)\|_{C(\overline{\Omega})}\leqslant C_1\|w_2(t)\|_{H^2(\Omega)}
\leqslant
C_2\|\Delta w_2(t)\|_{L_2(\Omega)}+
kC_2\|\Delta w_1(t)\|_{H^{\frac{1}{2}}(\Gamma_1)}+C_2\|w_2(t)\|_{L_2(\Omega)}
$$

$$
\leqslant C_2\|\Delta w_2(t)\|_{L_2(\Omega)}+
kC_3\|\Delta w_1(t)\|_{H^1(\Omega)}+C_2\|w_2(t)\|_{L_2(\Omega)}
$$

$$
\leqslant C_2\sqrt{\int\limits_{\Omega}\left\{(\Delta w_{1,x_1}(t))^2+(\Delta w_{1,x_2}(t))^2
+(\Delta w_2(t))^2\right\}dx}+
$$

$$
kC_3\sqrt{\int\limits_{\Omega}\left\{w^2_{2,x_1}(t)+w^2_{2,x_2}(t)+(\Delta w_1(t))^2\right\}dx}+
$$

$$
kC_3\sqrt{\int\limits_{\Omega}\left\{(\Delta w_{1,x_1}(t))^2+(\Delta w_{1,x_2}(t))^2+(\Delta w_2(t))^2\right\}dx}+
$$

$$
C_2\sqrt{\int\limits_{\Omega}\left\{w^2_{1,x_1}(t)+w^2_{1,x_2}(t)+w^2_{2}(t)\right\}dx}
$$

$$
\leqslant C_2\sqrt{\int\limits_{\Omega}\left\{(\Delta\varphi_{x_1})^2+(\Delta\varphi_{x_2})^2
+(\Delta\psi)^2\right\}dx}+
kC_3\sqrt{\int\limits_{\Omega}\left\{\psi^2_{x_1}+\psi^2_{x_2}+(\Delta\varphi)^2\right\}dx}+
$$

$$
kC_3\sqrt{\int\limits_{\Omega}\left\{(\Delta\varphi_{x_1})^2+(\Delta\varphi_{x_2})^2+(\Delta \psi)^2\right\}dx}+
C_2\sqrt{\int\limits_{\Omega}\left\{\varphi^2_{x_1}+\varphi^2_{x_2}+\psi^2\right\}dx}.
$$
Thus $\|w_2(t)\|_{C(\overline{\Omega})}$ is uniformly bounded for any $t$
because $k$ is near zero.

By virtue of (\ref{6})
initial conditions $(\varphi,\psi)$ will be the element of
$D(\mathfrak{A}^2)$ for any $k$.
If the coefficient $k$ is small enough, we achieve condition (\ref{4}).

\section{The second step of the control}

Now we have a task to drive the system to rest. A pair of functions
$$
(w|_{t=0}=w(T_1,x),w_t|_{t=0}=w_t(T_1,x))
$$
is considered to be new initial data for the problem (\ref{1})---(\ref{3}).
Bearing in mind that according to the fact proven above these initial conditions are sufficiently close to
$(C,0)$ in the norm of the space $\mathcal{C}^{4,3}(\overline{\Omega})$.
We shift now the solution $w$ (first step of the control) on the value $C$,
i.e. we change $w$ on $w+C$ and consider
the pair $(w|_{t=0}=w(T_1,x),w_t|_{t=0}=w_t(T_1,x))$ that is sufficiently close to
$(0,0)$.

Let us consider the domain $\Omega_{\delta}$, which is $\delta$-vicinity of the domain
$\Omega$. Also let take an arbitrary pair $(w_0(x),w_1(x))$ from the space
$\mathcal{C}^{4,3}(\overline{\Omega})$.
Consider an extension operator $E$. It is a linear continuous operator
from the space $\mathcal{C}^{4,3}(\overline{\Omega})$  to $\mathcal{C}^{4,3}(\overline{\Omega}_{\delta})$
such that the support of the extended pair $(w^e_0(x),w^e_1(x))$ and
its derivatives of 4th and 3th orders (respectively) inclusive belongs to
$\overline{\Omega}_{\delta}$. Moreover
$$
(w^e_0(x),w^e_1(x))=(w_0(x),w_1(x)),\:\:\mbox{if}\:\: x\in\overline{\Omega}.
$$
Note that, outside $\overline{\Omega}_{\delta}$, the
functions can be extended by zero to the whole plane. In
a more general case, $E$ was constructed in \cite{Madgenes}.

Extended in this way functions are denoted (as above) as
$w^e_0(x)$ and $w^e_1(x)$, according to D. L. Russell.

Let us consider the Cauchy problem for the equation of membrane's oscillations on a plane $R^2$:
\begin{equation}
\label{12}
w_{tt}(t,x)-\Delta w(t,x)=0,\quad (t,x)\in Q=(0,+\infty)\times R^2,
\end{equation}
\begin{equation}
\label{13}
w|_{t=0}=w^e_0(x),\quad w_t|_{t=0}=w^e_1(x),\quad x\in R^2.
\end{equation}

It is known that the solution to the problem (\ref{12}), (\ref{13}) has the form (Poisson's formula):
\begin{equation}
\label{14}
w(t,x)=\frac{\partial}{\partial t}\left(\frac{1}{2\pi}\int\limits_{|y-x|
< t}\frac{w^e_0(y)dy}{\sqrt{t^2-|y-x|^2}}\right)+
\frac{1}{2\pi}\int\limits_{|y-x|
< t}\frac{w^e_1(y)dy}{\sqrt{t^2-|y-x|^2}}.
\end{equation}

We use the formula (\ref{14}) for estimating the absolute value of the solution $w(t,x)$ uniformly by
the initial data. The absolute value of $w(t,x)$ is estimated in case $x\in\overline{\Omega}_{\delta}$.
We compute the first derivative with respect to $t$ in the right part of (\ref{14}):
\begin{equation}
\label{14.1}
w(t,x)=\frac{1}{2\pi t}\int\limits_{|y-x|
< t}\frac{w^e_0(y)+(y-x)\cdot\nabla w^e_0(y)}{\sqrt{t^2-|y-x|^2}}dy+
\frac{1}{2\pi}\int\limits_{|y-x|
< t}\frac{w^e_1(y)dy}{\sqrt{t^2-|y-x|^2}}.
\end{equation}
As initial data $(w^e_0(x),w^e_1(x))$ have a compact support then there is large enough time $t^{\ast}>0$
such that for any $t>t^{\ast}$ and for any $x\in\overline{\Omega}_{\delta}$ we obtain
\begin{equation}
\label{14.2}
w(t,x)=\frac{1}{2\pi t}\int\limits_{\Omega_{\delta}}
\frac{w^e_0(y)+(y-x)\cdot\nabla w^e_0(y)}{\sqrt{t^2-|y-x|^2}}dy+
\frac{1}{2\pi}\int\limits_{\Omega_{\delta}}
\frac{w^e_1(y)dy}{\sqrt{t^2-|y-x|^2}}.
\end{equation}
Note that we choose $t$ such as $t^2-|y-x|^2\geqslant\alpha>0$ for any
$x,y\in\overline{\Omega}_{\delta}$.

The following rough evaluation follows from the explicit form of (\ref{14.2}):
\begin{equation}
\label{15}
\|w(t,\cdot)\|_{C^4(\overline{\Omega}_{\delta})}\leqslant
\frac{C_1}{t}\|w^e_0\|_{C^4(R^2)}+\frac{C_2}{t}\|w^e_1\|_{C^3(R^2)}.
\end{equation}

Differentiating $w(t,x)$ with respect to $t$, we obtain the rough estimate in the space of the pair of functions
$\mathcal{C}^{4,3}(\overline{\Omega}_{\delta})=C^4(\overline{\Omega}_{\delta})
\times C^3(\overline{\Omega}_{\delta})$
\begin{equation}
\label{16}
\|(w(t,\cdot),w_t(t,\cdot))\|_{\mathcal{C}^{4,3}(\overline{\Omega}_{\delta})}
\leqslant\frac{M}{t}\|(w^e_0,w^e_1)\|_{\mathcal{C}^{4,3}(R^2)},\quad t>t^{\ast},
\end{equation}
where a number $M$ does not depend on initial data.

Further we use the method described in \cite{Russell} and applied to problems of the boundary controllability
for a wave equation.

Let us consider some initial conditions $w_0(x)$ and $w_1(x)$, $x\in\Omega$.
We extend them to $R^2$ by means of a linear bounded operator $E$.
Then we obtain $(w^e_0,w^e_1)=E(w_0,w_1)$. And the Cauchy problem $(\ref{12})$, $(\ref{13})$ arises.
Let $w^s(t,x)$ be the solution to this Cauchy problem. Now consider any large enough time $t=T_2$. We get $(w^s(T_2,x),w^s_t(T_2,x))\in C^{4,3}(\overline{\Omega})$. The restriction of the function $w^s(T_2,x)$ and its derivative on the domain $\Omega$ should be considered. It is obvious that in virtue of (\ref{16}) the following estimate is correct for $t=T_2$
\begin{equation}
\label{17}
\|(w^s(T_2,\cdot),w^s_t(T_2,\cdot))\|_{\mathcal{C}^{4,3}(\overline{\Omega})}
\leqslant\frac{M}{T_2}\|(w^e_0,w^e_1)\|_{\mathcal{C}^{4,3}(R^2)}.
\end{equation}
Let by definition
$(w^{s,e}_0(T_2,x),w_1^{s,e}(T_2,x))=E\left(w^s(T_2,x)|_{\Omega},w^s_t(T_2,x)|_{\Omega}\right)$.
Now let us have a look at the inverse Cauchy problem with initial conditions
\begin{equation}
\label{18}
w(t,x)|_{t=T_2}=-w^{s,e}_0(T_2,x)\quad w_t(t,x)|_{t=T_2}=-w_1^{s,e}(T_2,x).
\end{equation}
Let $w^i(t,x)$ be the solution to the inverse Cauchy problem with conditions (\ref{18}).
In virtue of invertibility of the equation (\ref{1}) with respect to $t$ the following estimate takes place:
\begin{equation}
\label{19}
\|(w^i(0,\cdot),w^i_t(0,\cdot))\|_{\mathcal{C}^{4,3}(\overline{\Omega})}
\leqslant\frac{M}{T_2}\|(w^{s,e}_0(T_2,x),w_1^{s,e}(T_2,x))\|_{\mathcal{C}^{4,3}(R^2)}.
\end{equation}

Obviously the solution of the Cauchy problem with initial conditions such as
\begin{equation}
\label{20}
w|_{t=0}=w_0^e(x)+w^i(0,x),\quad w_t|_{t=0}=w_1^e(x)+w^i_t(0,x),\quad x\in R^2,
\end{equation}
identically equals zero in $\Omega$ as well as its first derivative with respect to $t$ at the time
$t=T_2$. Now let us consider the restriction of the right parts of (\ref{20}) in the domain $\Omega$.
We regard the initial conditions (the restriction of right parts of (\ref{20}) in the domain $\Omega$) in the problem of boundary controllability:
\begin{equation}
\label{21}
w|_{t=0}=w_0(x)+w^{i,r}(0,x),\quad w_t|_{t=0}=w_1(x)+w^{i,r}_t(0,x),\quad x\in\Omega.
\end{equation}
Note that it is the value of the corresponding solution to the Cauchy problem in $R^2$ with the
initial conditions (\ref{20}) to determine the required control function on the boundary of $\Omega$.

A pair $(w^{i,r}(0,x),w^{i,r}_t(0,x))$ is derived from pair $(w_0(x),w_1(x))$  by means of applying a linear continuous operator, let us denote it as $L$, with the norm less than $1$
(consequence from estimates (\ref{17}) and (\ref{19})).
Obviously the sums in right parts (\ref{21}) generate all elements of the space
$\mathcal{C}^{4,3}(\overline{\Omega})$. Indeed, (\ref{21}) can be written as:
\begin{equation}
\label{22}
(I+L)(w_0(x),w_1(x))=(w|_{t=0},w_t|_{t=0}),
\end{equation}
where $I$ is the identical operator. Hence, as $\|L\|<1$, so the operator
$I+L$, which acts from $\mathcal{C}^{4,3}(\overline{\Omega})$ to itself,
is invertible.

Now let us represent the control function (second step) in the following form:
$$
u^{(2)}(t,x)=\frac{\partial}{\partial\nu}PK^t_+\left[\left(I+(-K^{T_2}_-)ERK^{T_2}_+\right)
E\left(I+R(-K^{T_2}_-)ERK^{T_2}_+E\right)^{-1}(w|_{t=0},w_t|_{t=0})\right],\:
x\in\partial\Omega,
$$
where $R$ is a restriction from $R^2$ to $\Omega$ and $K^{T_2}_+$, $K^{T_2}_-$
are resolving operators of the Cauchy problem and $P$ is a projection:
$(a,b)\mapsto a$. We write minus before $K^{T_2}_-$
because of (\ref{18}).

Thus we have proven that the system with smooth initial conditions can be driven to rest by means of extending them on the full plane. It is the method to extend which determines a program of the boundary control. Let us show now that if the initial conditions have small enough absolute values, we can drive the system to rest by means of a boundary control which has a small absolute value.

We regard that in the problem (\ref{1})---(\ref{3}) the value of the solution $w(t,x)$ and the value of its
derivative $w_t(t,x)$
at $t=T_1$ are small enough in norms of spaces $C^4(\overline{\Omega})$  and
$C^3(\overline{\Omega})$ respectively.

Let $(w|_{t=0},w_t|_{t=0})$ be rewritten according to the formula
(\ref{21}). As continuous operator $I+L$ invertible, so according to Banach's theorem
an invertible operator is continuous too. Hence choosing $(w|_{t=0},w_t|_{t=0})$ sufficiently small,
we can make $(w_0(x),w_1(x))$ be sufficiently small as well. Now let consider the sums (\ref{21}), which determine data $(w|_{t=0},w_t|_{t=0})$. Extending these sums on the whole plane by the method above, we obtain initial data $(\ref{20})$.

Bearing in mind that supports of functions $w^e_0(x)$ and $w^e_1(x)$ are in $\overline{\Omega}_{\delta}$,
and supports of their derivatives with respect to all variables (including the third
and the second orders respectively) are located in $\overline{\Omega}_{\delta}$ too.
The solution $w^s(t,x)$ has a compact support which is located in some bounded domain $G_t$ in $R^2$ at each moment $t$ because of the finite speed of the wave propagation. Let us take a sufficiently large
circle $D$ such as $\overline{G}_t\subset D$, $t\in[0,T_2]$.
In this case function $w^s(t,x)$ is thought as a solution of initial boundary value problem at the domain $D$ with the homogeneous Dirichlet condition for $t\in[0,T_2]$. In virtue of the corresponding smoothness of initial conditions we obtain:
$w^s(t,x)\in C([0,T_2];H^4(D))$ and $w^s_t(t,x)\in C([0,T_2];H^3(D))$.
Then the energy conservation law takes place:
$$
\int\limits_{D}\left\{\left(w^s_{x_1}(t,x)\right)^2+\left(w^s_{x_2}(t,x)\right)^2+
\left(w^s_t(t,x)\right)^2\right\}dx=
$$
\begin{equation}
\label{25}
\int\limits_{D}\left\{\left(\frac{\partial w^e_0(x)}{\partial x_1}\right)^2+
\left(\frac{\partial w^e_0(x)}{\partial x_2}\right)^2+\left(w^e_1(x)\right)^2\right\}dx,
\quad t\in[0,T_2].
\end{equation}

Now differentiating the equation (\ref{12})
and initial conditions (\ref{13}) with respect to variables $x_1$, $x_2$,  we obtain the estimate
$$
\|w^s(t,\cdot)\|^{\prime}_{H^3(D)}\leqslant
\|w^e_0\|_{H^3(\Omega_{\delta})}+\|w^e_1\|_{H^2(\Omega_{\delta})},
$$
where $\|\cdot\|^{\prime}_{H^3(D)}$ is a seminorm (term
$$
\int\limits_{D}(w^s(t,x))^2dx
$$
is absent).
The last statement is true because derivatives (of the second order in this case) of function $w^s(t,x)$
are identically zero at domain $D\setminus\overline{G}_t$ and hence they are solutions of differentiated initial boundary value problem with the homogeneous Dirichlet condition at the boundary of the domain $D$.

Then the seminorm
$\|\cdot\|^{\prime}_{H^3(D)}$ is a norm. Therefore we obtain
$$
\|w^s(t,\cdot)\|_{H^3(D)}\leqslant
C_F\|w^e_0\|_{H^3(\Omega_{\delta})}+C_F\|w^e_1\|_{H^2(\Omega_{\delta})}.
$$

Taking into account the last estimate and the Sobolev embedding theorem we get
\begin{equation}
\label{26}
\|w^s(t,\cdot)\|_{C^1(\overline{\Omega})}\leqslant
C_S\|w^e_0\|_{H^3(\Omega_{\delta})}+C_S\|w^e_1\|_{H^2(\Omega_{\delta})}.
\end{equation}

Summing up, it is proven that the solution $w^s(t,x)$ can be made sufficiently small in the norm $C^1(\overline{\Omega})$ for any $t\in[0,T_2]$.
The same argument may be applied to the solution of the inverse Cauchy problem with initial conditions $-w^{s,e}_0(T_2,x)$ and $-w^{s,e}_1(T_2,x)$.
In this case it is important that functions $w^{s,e}_0(T_2,x)$ and $w^{s,e}_1(T_2,x)$
in virtue of inequality (\ref{16}) are "small" in $\mathcal{C}^{4,3}$, if $w^e_0(x)$ and $w^e_1(x)$
are "small". Hence the restriction of the normal derivative of the solution to the Cauchy problem (\ref{12}),
(\ref{13}) on the boundary of $\Omega$ (Neumann condition of the problem of controllability) is less than given
$\varepsilon$ with respect to absolute value.
The latter means that the required restriction (\ref{4}) on the control function $u(t,x)$ is satisfied.

\end{document}